\documentclass[11pt, oneside]{article}   	
\usepackage{geometry}                		
\newtheorem{thm}{Theorem}

\newtheorem{conj}{Conjecture}

\usepackage{graphicx}				
\usepackage{amssymb}

\title{Ideals of generic forms}
\author{Ralf Fr\"oberg}
\date{}
\begin{document}
\maketitle
\begin{abstract}
We determine the Hilbert series of some classes of ideals generated by generic forms of degree two and three, and investigate the
difference to the Hilbert series of ideals generated by powers of linear generic forms of the corresponding degrees. 
\end{abstract}

Adress\footnote
{Department of mathematics, Stockholm university, S-10691 Stockholm Sweden\\
ralf.froberg@math.su.se}
\section{Introduction}

Let $R=k[x_1,\ldots,x_n]/(f_1,\ldots,f_r)$, where $k$ is a field and $f_i$ is a homogeneous polynomial of degree $d_i$, $i=1,\ldots,r$. 
Then $R$ is a graded algebra, $R=\oplus_{i\ge0}R_i$, where $R_i$ are finite dimensional vector spaces over $k$.
The generating series
of these dimensions, $H_R(t)=\sum_{i\ge0}\dim_kR_it^i$ is called the Hilbert series of $R$. 
It is shown in \cite{Fr-Lo} that for fixed $n,d_1,\ldots,d_r$ there is only a finite number of Hilbert series. If we think of $f_i$ as a linear 
combination of all monomials of degree $d_i$, we can view the ideal as a point in ${\Bbb A}^N_k$, where $N=\sum_{i=1}^r {n+d_i-1\choose d_i}$.
It is further shown in \cite{Fr-Lo} that if $k$ is infinite, there is an nonempty Zariski open set in ${\Bbb A}^N_k$ where the Hilbert series is constant,
so "almost all" ideals give the same Hilbert series. The series is the smallest among algebras of this form.
We call the ideals in this set {\it generic}. If $k={\Bbb C}$, we can think of the case when all $N$ coefficients are algebraically independent over ${\Bbb Q}$.

There is a longstanding conjecture on what this series is.

\begin{conj}\label{conj}\cite{Fr} The Hilbert series of $k[x_i,\ldots,x_n]/(f_1,\ldots,f_r)$, where $k$ is infinite and $f_i$ is a generic form of degree $d_i$
is $[\frac{\prod_{i=1}^r(1-t^{d_i})}{(1-t)^n}]$. Here $[\sum_{i=0}^\infty a_it^i]=\sum_{i=0}^\infty b_it^i$, where
$b_i=a_i$ if $a_j>0$ for all $j\le i$ and $b_i=0$ otherwise.
\end{conj}
The conjecture is trivially true for $r\le n$. It is proved for $n\le3$, \cite{An}, for $r=n+1$, \cite{St}, and for lots of other special cases, see \cite{Fr-Lu}.
From now on we suppose that $d_i=d$ for all $i$, so the conjectured series is $[\frac{(1-t^d)^r}{(1-t)^n}]$.
This note is an amplification of \cite{Fr-Ho}, where the conjecture was proved for $d=2,n\le11$ and $d=3,n\le8$. In \cite{Fr-Ho} was also
investigated when the series for $r$ generic forms of degree $d$ equals the series for $r$ $d$'th powers of generic linear forms.

\section{Results}
Now we extend the results from \cite{Fr-Ho}.
\begin{thm} Suppose ${\rm char} k=0$. Then conjecture \ref{conj} is true if if $(d,n)=(2,12),(2,13),(3,9)$.
\end{thm}
{\bf Proof} It is known that the true Hilbert series is the smallest possible among algebras of these respective types, and
greater or equal to the conjectured lexicographically, see \cite{Fr}. Thus it is sufficient to find one example for the 
parameters $(n,r,d)$ which gives equality,
to prove the theorem for those parameters. This has been done with computer for the ideals $(x_1^d,\ldots,x_n^d)+J$, where $J$ is generated by $r-n$ random
forms of degree $d$. The calculations where made in Macaulay2, \cite{M2}.

\bigskip\noindent
More interesting is perhaps next theorem. Let the Hilbert series of $k[x_1,\ldots,x_n]/(l_1^d,\ldots,l_r^d)$,  where $l_i$ are generic
linear forms, be $Q_{n,r,d}(t)$, and let $F_{n,r,d}(t)=[\frac{1-t^d)^r}{(1-t)^n}]$ , which in the cases we study is the correct Hilbert series 
for $k[x_i,\ldots,x_n]/(f_1,\ldots,f_r)$,
where the $f_i$'s are generic forms of degree $d$.
\begin{thm}
For $n=12$, $Q_{n,r,2}(t)=F_{n,r,2}(t)$ unless $r=n+2$. For $n=13$, $Q_{n,r,2}(t)=F_{n,r,2}(t)$ unless $r=n+2$ or $r=n+3$.
For $n=9$, $Q_{n,r,3}(t)=F_{n,r,3}(t)$ unless $r=n+2$ or $r=n+3$.
\end{thm}
{\bf Proof} Now we use the ideals $(x_1^d,\ldots,x_n^d,f_1^d,\ldots,f_{r-n}^d)$, where the $f_i$'s are random linear forms. They all have the conjectured series.

\medskip
The following claim is checked several times on computer. That $Q(t)\ne F(t)$ in our calculations might depend on bad luck (which is highly unlikely
especially since we made the calculations with different random linear forms with the same answer). Thus, it is not a theorem, it is a conjecture.
\begin{conj} $Q_{12,14,2}(t)-F_{12,14,2}(t)=64t^6$, $Q_{13,15,2}(t)-F_{13,15,2}(t)=13t^6+t^7$, $Q_{13,16,2}(t)-F_{13,16,2}(t)=t^6$, $Q_{9,11,3}(t)-F_{9,11,3}(t)=t^8+154t^9+t^{10}$,
$Q_{9,12,3}(t)-F_{9,12,3}(t)=12t^8$.
\end{conj}
For the geometric implications of $Q_{n,r,d}(t)=F_{n,r,d}(t)$, and a conjecture when it is true, see \cite{Ia}, \cite{Ch}. 

\medskip
{\bf Declarations of interest}: None

\end{document}